\documentclass[leqno, letterpaper]{article}
\usepackage{amssymb, amsmath, amsthm, amscd, mathrsfs}
\usepackage{leqno, tocbibind}

\usepackage[pdftex]{hyperref}
\hypersetup{
  colorlinks   = true, 
  urlcolor     = black, 
  linkcolor    = black, 
  citecolor   =  black
}

\newcommand{\on}{\operatorname}

\newcommand{\mf}{\mathfrak}
\newcommand{\mb}{\mathbf}
\newcommand{\br}{\buildrel}

\newcommand{\be}{\begin{equation}}
\newcommand{\ee}{\end{equation}}

\newcommand{\CC}{{\mathbb C}}

\newcommand{\PP}{{\mathbb P}}
\newcommand{\QQ}{{\mathbb Q}}

\newcommand{\ZZ}{{\mathbb Z}} 
\newcommand{\mbfk}{{\mathbf k}}

\newcommand{\Cc}{{\mathcal C}} 
 
\newcommand{\Ec}{{\mathcal E}} 
\newcommand{\Fc}{{\mathcal F}} 
\newcommand{\Gc}{{\mathcal G}} 
\newcommand{\Hc}{{\mathcal H}} 
\newcommand{\Ic}{{\mathcal I}} 
 
\newcommand{\Kc}{{\mathcal K}}

\newcommand{\Oc}{{\mathcal O}} 
 
\newcommand{\Qc}{{\mathcal Q}}

\newcommand{\Vc}{{\mathcal V}}

\newcommand{\Ms}{{\mathscr M}}

\newcommand{\Hom}{\on{Hom}}

\newcommand{\Id}{\on{Id}}

\newcommand{\Spec}{\on{Spec}}
\newcommand{\Grass}{\on{Grass}}
\newcommand{\Quot}{\on{Quot}}

\newcommand{\Ind}{\on{Ind}}

\newcommand{\QCoh}{\mathcal{QC}oh}
\newcommand{\Coh}{\mathcal{C}oh}

\newcommand{\Proj}{\on{Proj}}
\newcommand{\Schb}{{\mb{Sch}}}
\newcommand{\Setb}{{\mb{Set}}}
\newcommand{\Grqc}{\underleftarrow{\on{G}}}
\newcommand{\Qqc}{\underleftarrow{\on{Q}}}

\newcommand{\Rat}{\on{Rat}}

\newcommand{\ind}{\varinjlim}
\newcommand{\pro}{\varprojlim}

\newcommand{\lra}{\longrightarrow}

\newcommand{\xrightarrowdbl}[2][]{%
  \xrightarrow[#1]{#2}\mathrel{\mkern-14mu}\rightarrow
}
\newtheorem{thm}[equation]{Theorem}

\newtheorem{lem}[equation]{Lemma}
\newtheorem{prop}[equation]{Proposition}

\newtheoremstyle{example}{\topsep}{\topsep}%
     {}
     {}
     {\bfseries}
     {.}
     {2pt}
     {\thmname{#1}\thmnumber{ #2}\thmnote{ #3}}

   \theoremstyle{example}
   \newtheorem{nota}[equation]{Notation}
   \newtheorem{Defi}[equation]{Definition}
   \newtheorem{rem}[equation]{Remark}

   \newtheorem{exa}[equation]{Example}
   
\title{The \emph{quot} Functor of a Quasi-coherent Sheaf}
\author{Gennaro Di Brino}
\date{}
\begin{document}

\maketitle
\pagenumbering{arabic}

\begin{abstract}
We build an infinite dimensional scheme parametrizing isomorphism classes of {\em coherent} quotients of a \emph{quasi-coherent} sheaf on a projective scheme. The main tool to achieve the construction is a version of Grothendieck's Grassmannian embedding combined with a result of Deligne, realizing quasi-coherent sheaves as ind-objects in the category of quasi-coherent sheaves of finite presentation. We end our treatment with the discussion of a special case in which we can retain an analog of the Grassmannian embedding.
\end{abstract}

\tableofcontents

\vfill\eject

\section{Introduction}

\paragraph*{The original construction and our leading question.} Since their introduction in \cite{Gr}, \emph{Quot} schemes have played a fundamental role in algebraic geometry and, in particular, in deformation theory. For instance, they provide natural compactifications of spaces of morphisms between certain schemes (\cite{Gr}), they are used to give a presentation of the stack of coherent sheaves over a projective scheme (\cite{LMB}), and their derived version (\cite{CFK}, \cite{Pri}) is of fundamental importance in derived algebraic geometry.\\
Recall that, given a Hilbert polynomial $h$, a projective scheme $X$ over an algebraically closed field $\mbfk$ and a quasi-coherent sheaf $\Ec$, one defines the contravariant functor $quot^X_{h}(\Ec_{(-)}):(\Schb/\mbfk)^{\circ}\rightarrow\Setb$ as

\begin{equation}\label{e:qf}
quot^X_{h}(\Ec_T):=\left\{\Kc\subset\Ec_T\,\, | \begin{array}{cc}\Ec_T/\Kc\,\, {is}\,\,{coherent,}\,\,{flat}\,\,{over}\,\,\Oc_T,\\{and}\,\,{has}\,\,{Hilbert}\,\,{polynomial}\,\, h
\end{array}\right\}
\end{equation}

where $\Ec_T:=\Ec\otimes_{\mbfk}\Oc_T$, and as pullback on morphisms. Note also that coherence, flatness and the Hilbert polynomial condition are imposed on \emph{quotients}. In his original outline of the construction, Grothendieck proves the representability of the above functor only in the case when $\Ec$ is \emph{coherent}. \\
The question we address in this work concerns the representability of the quot functor when the sheaf $\Ec$ in \eqref{e:qf} is only assumed to be quasi-coherent. In particular, our main result (Theorem \ref{T:iquot-rep}) implies the statement below.

\begin{thm}
Let $\Ec$ be a quasi-coherent $\Oc_X$-module and let $X$ and $h$ be as above. Then there is a scheme $\Quot_h^X(\Ec)$ representing the functor $quot^X_{h}(\Ec_{(-)})$.
\end{thm}

In the following two paragraphs we give an idea of our construction of the scheme $\Quot_h^X(\Ec)$ and outline the other results contained in this work. In the rest of this Introduction $\Ec$ will always denote a quasi-coherent sheaf over the projective scheme $X$ and we may sometimes refer to $\Quot_h^X(\Ec)$ as the \emph{quasi-coherent Quot scheme}.

\paragraph*{The filtering schematic Grassmannian.} The main idea in Grothendieck's paper is that the representability of the \emph{quot} functor and the resulting universal property of the \emph{Quot} scheme are inherited from the corresponding properties of a certain Grassmannian, in which the \emph{quot} functor lives. Motivated by this, in Section \ref{S:grass} we provide a filtering construction for the Grassmannian. More in detail, if $\Grass_n(\Fc)$ denotes the Grassmannian of locally free rank $n$ quotients of a quasi-coherent sheaf $\Fc$ \cite{EGAS}, we prove the following proposition (see Lemma \ref{L:grsub}).

\begin{prop}\label{P:grass}
Let $\Fc$ be a quasi-coherent sheaf over a scheme $S$. Then $\Grass_n(\Fc)$ is the filtering inductive limit over $i$ of an increasing sequence of quasi-compact open subschemes $\left(\Grqc(n,\Fc)_i\right)_{i\in I}$.
\end{prop}

The schemes $\Grqc(n,\Fc)_i$ in the statement are constructed as projective limits of diagrams consisting of certain subschemes of Grassmannians of finite type and affine morphisms between them (see Lemma \ref{L:affmorphs} for the proof of affineness).
Also, note that a crucial ingredient in the proof of Proposition \ref{P:grass} is the following theorem of Deligne.

\begin{lem}{\cite{Del}}
Let $X$ be a quasi-compact quasi-separated scheme (not necessarily Noetherian). Then the category $\QCoh(X)$ of quasi-coherent sheaves on $X$ is equivalent to that of ind-objects in the category of quasi-coherent sheaves of finite presentation on $X$.
\end{lem}

For the sake of the reader, we briefly review the concepts involved the above statement in Subsection \ref{SS:indobj}.

\paragraph*{Two uses of the Grassmannian embedding.} The filtering construction from Section \ref{S:grass} will be of twofold interest to us. First, in building the $Quot$ scheme of $\Ec$ we will proceed in a way that follows the same ``ind-pro principle'', to get an infinite dimensional analog of the classical scheme of \cite{Gr} which we denote again by $\Quot_{h}^X(\Ec)$. More precisely, we construct a candidate for the scheme representing the functor $quot^X_{h}(\Ec_{(-)})$ as a filtering inductive limit of certain schemes denoted $\Qqc(h,\Ec)_i$. In order to obtain the $\Qqc(h,\Ec)_i$'s, we take the projective limit of a filtering projective system consisting of some open subschemes of ordinary \emph{Quot} schemes and affine morphisms between them. Roughly speaking, the affineness of such morphisms will be proved by viewing them as restrictions of morphisms between Grassmannians (Lemma \ref{L:qaffmorph}).\\
In the last part of this work, we introduce \emph{uniformly regular} sheaves over a projective scheme $X$. These are quasi-coherent sheaves all of whose coherent approximations in the sense of Proposition \ref{P:Deligne} (and the following Remarks) have Castelnuovo-Mumford regularities \cite{Mum} bounded by a given integer $m$. This said, the other way in which we use the filtering construction from Proposition \ref{P:grass} is to show that $\Quot_{h}^X(\Ec)$ can be embedded in some schematic Grassmannian. The precise statement is as follows (see Definition \ref{D:quasiclosed}).

\begin{thm}\label{P:urgrassemb}
Let $\Ec$ be a uniformly $m$-regular quasi-coherent sheaf on a projective $\mbfk$-scheme $X$. Then there is a quasi-closed embedding 
$$
\Quot^X_h(\Ec)\hookrightarrow \Grass_{h(m)}\left(\ind_j H^0(X,\Ec^j(m))\right).
$$
\end{thm}

\paragraph*{Further directions of research.} Assume $\mbfk=\CC$. Recall that given nonnegative integers $d, r$ and $m$, with $r<m$, the ordinary \emph{Quot} scheme $\Quot^{\PP^1}_{r,d}(\Oc_{\PP^1}^m)$ of quotients of $\Oc_{\PP^1}^m$ of rank $r$ and degree $d$, can be used to compactify the space 
$$
\Ms_d^m:=\Rat_d(\PP^1,\Grass_r(\Oc_{\PP^1}^m)),
$$
of maps of degree $d$ from $\PP^1$ to $\Grass_r(\Oc_{\PP^1}^m)$. To see this, one can note that giving a morphism $\PP^1\rightarrow\Grass_r(\Oc_{\PP^1}^m)$ of degree $d$ is equivalent to the datum of a quotient bundle of $\Oc_{\PP^1}^m$ of rank $r$ and degree $d$. Similarly, letting $m$ approach $\infty$, one could use the \emph{Quot} scheme we construct here to study the topology of the spaces 
$$
\Ms_d^{\infty}= \Rat_d(\PP^1,\Grass_r(\Oc_{\PP^1}^{\oplus\infty})),
$$
of rational curves of degree $d$ in the schematic Grassmannian.\

\paragraph*{Relations to other work.} An instance of the quasi-coherent \emph{Quot} scheme of length $1$ quotients appeared in the work of S. Kleiman \cite{Kle}, we refer the reader to Remark \ref{R:Kleiman} for further details.\\
More recently, keeping $\Ec$ quasi-coherent, a generalization of Grothendieck's theorem in another direction has been carried out by R. Skjelnes, applying different techniques from the ones used here. More precisely, in \cite{Skj} it is showed that if one replaces the projective scheme $X$ with an \emph{algebraic space} and requires the quotients to be flat, with \emph{finite support} and of \emph{relative rank} $n$, then the object representing the quot functor is again an algebraic space (\cite[Theorem 6.11]{Skj}).\\
Finally, while a second version of this paper was being prepared, the author learned from David Rydh that a construction for general \emph{quot} functors was being carried out independently by Hall and Rydh in what is now available as \cite{HR}. In their case, starting from an algebraic stack $X$ with finite diagonal and a quasi-coherent sheaf on it, one obtains a separated algebraic space representing the \emph{quot} functor in the sense of algebraic spaces.

\paragraph*{Acknowledgments.}

The author would like to thank his Ph.D. advisor, Mikhail Kapranov, for suggesting the problem and for his constant support and encouragement. Besides, the author wishes to express his gratitude to David Rydh for kindly sharing a draft of his joint work with Hall, and to Roy Skjelnes and Angelo Vistoli for helpful comments and suggestions. Finally, the author is grateful to the Mathematics Department of the KTH, Stockholm, for the warm hospitality during the preparation of a second version of this paper.

\numberwithin{equation}{subsection}

\section{Some background material}

We recall here some results that we will need in the rest of our work.

\subsection{Limits and quasi-compact schemes}

Recall that a \emph{filtering inductive limit} is defined as the inductive limit of a covariant functor $F:\Ic\rightarrow \Cc$ where $\Cc$ is any category and $\Ic$ is a filtering poset, i. e., every two objects have a common upper bound. In case no confusion can arise, we may refer to $\Ic$ simply as the indexing category. A \emph{filtering projective limit} is defined dually, assuming that the functor $F$ above is contravariant.
\vskip 3mm
The following result will be used to construct quasi-compact schemes that will form two filtering covers for the schematic Grassmannian and the quasi-coherent Quot scheme, respectively. The statements from which it follows are scattered throughout \cite{EGAIV3} to which the reader is referred. We remark that affineness of the morphisms in the statement below ensures that the limit is a scheme, whereas quasi-compactness of the limit follows from Tychonoff's Theorem.\\

\begin{prop}\label{P:qcam}
Let $S$ be a scheme and let $(X_i)_i$ be a filtering projective system of quasi-compact $S$-schemes and affine morphisms between them. Then the projective limit
$$
X:=\pro_{i}\,^{\Schb/S}X_i
$$
is a quasi-compact $S$-scheme.
\end{prop}

\subsection{Ind-objects and a theorem of Deligne}\label{SS:indobj}

An \emph{ind-object} in a category $\Cc$ is by definition a filtering inductive limit in $\widehat{\Cc}:=\Fc un(\Cc^{\circ},\Setb)$ of presheaves of sets of $\Cc$. One denotes by $\Ind(\Cc)$ the full subcategory of $\widehat{\Cc}$ whose objects are the ind-objects of $\Cc$ (see \cite{Del} or more extensively \cite{SGAIV2} for further details).

The proposition below will be of crucial importance in what follows.

\begin{prop}{\cite[Prop. 2]{Del}}\label{P:Deligne}
Let $X$ be a quasi-compact quasi-separated scheme (not necessarily Noetherian). Then the category $\QCoh(X)$ of quasi-coherent sheaves on $X$ is equivalent to that of ind-objects in the category of quasi-coherent sheaves of finite presentation on $X$.
\end{prop}
 
\begin{rem}\label{R:univincl}
A quasi-coherent sheaf $\Ec$ on a quasi-compact quasi-separated scheme $X$ is therefore given by the inductive limit in $\QCoh(X)$ of a filtering inductive diagram 
\be\label{e:diagram}
\left(\Ec^{i},\alpha^{i,j}:\Ec^{i}\rightarrow \Ec^{j}\right)_{i,j\in\Ic} 
\ee
of finitely-presented sheaves and morphisms between them.
\end{rem}

Next, the lemma below shows that on a Noetherian scheme the notion of sheaf of finite presentation reduces to that of coherent sheaf.

\begin{lem}\label{L:fpqconnoethiscoh}
Let $X$ be a Noetherian scheme and let $\Fc$ be a finitely presented quasi-coherent sheaf on $X$. Then $\Fc$ is coherent.
\end{lem}

\noindent {\sl Proof:} Using the definition of coherent sheaves from \cite{EGAS}, we need to check that $\Fc$ is of finite type and that for all opens $U$ and for all integers $n$, the kernel of any homomorphisms $\Oc_X^n|U\rightarrow \Fc|U$ is of finite type. The first condition being automatically satisfied, we need only check that kernels are of finite type. But this follows from the fact that the submodules of a Noetherian module are finitely generated, which in turn holds since a finitely generated module over a Noetherian ring, such as $\Oc_X|U$, is again Noetherian.
\qed

\begin{rem}\label{R:qcoh-equivalent-to-indcoh-on-noetherian-schemes}
Assuming that $X$ is Noetherian, Proposition \ref{P:Deligne} together with Lemma \ref{L:fpqconnoethiscoh} imply that there is an equivalence of categories between the category $\QCoh(X)$ and $\Ind(\Coh(X))$, the ind-category of the category of coherent sheaves on $X$. We will sometimes refer to the coherent sheaves $\Ec^i$ in the diagram \eqref{e:diagram} as to the $i$-th \emph{coherent approximation} of $\Ec$.
\end{rem}

\subsection{Representable functors}

Let $S$ be a scheme, for any contravariant functor 
$$
\varphi: (\Schb/S)^{\circ}\rightarrow \Setb,
$$
the datum $U\mapsto \varphi(U)$, with $U\subset X$ open, defines in the obvious way a presheaf of sets on every $S$-scheme $X$. The functor $\varphi$ is then called a sheaf in the Zariski topology if for any $S$-scheme $X$ and for any cover $\{U_{\alpha}\}_{\alpha}$ of $X$ the sequence
$$
\varphi(X)\rightarrow \prod_{\alpha} \varphi(U_{\alpha})\rightrightarrows\prod_{\alpha, \beta}\varphi(U_{\alpha}\times_X U_{\beta})
$$
is exact, i.e., if $\{f_{\alpha}\in\varphi(U_{\alpha})\}_{\alpha}$ are such that $f_{\alpha}$ and $f_{\beta}$ agree in $\varphi(U_{\alpha}\times_X U_{\beta})$ then there exists a unique $f\in \varphi(X)$ mapping to each $f_{\alpha}$ via the restriction $\varphi(X)\rightarrow \varphi(U_{\alpha})$.
\vskip 3mm
It is well known that being a sheaf in the Zariski topology is a necessary condition for the functor $\varphi$ to be representable. Moreover, once one is able to prove that $\varphi$ is a sheaf one can reduce to showing representability over the category of affine schemes. The following Lemma provides a criterion for $\varphi$ to be representable, we recall it here for the reader's convenience (see, e.g., \cite[Appendix E]{Ser}).

\begin{lem}
Let $S$ and $\varphi$ be as above. Assume that
\begin{itemize}
\item[a)] $\varphi$ is a sheaf;
\item[b)] $\varphi$ admits a covering by representable open subfunctors $\varphi_{\alpha}$.
\end{itemize}
Then $\varphi$ is representable by an $S$-scheme.
\end{lem}

\subsection{Preliminaries to the classical \emph{Quot} scheme}\label{SS:prelquot}

In this subsection we summarize Grothendieck's construction of the \emph{Quot} scheme and the main results needed. Besides \cite{Gr}, the more extensive treatments we refer the reader to are \cite{Mum}, \cite{Vie} and \cite{HL}. The first section of \cite{CFK} also contains a brief outline of the construction, and some of the statements we will need are closer in spirit to those.\\
Let $\mbfk$ be an algebraically closed field and let $X$ be a projective $\mbfk$-scheme, considered together with a fixed very ample invertible sheaf $\Oc_X(1)$.

\vskip 3mm

The theorem of Serre (see \cite{FAC} or \cite{CFK}) below already contains the notion of what is usually referred to as \emph{Castelnuovo-Mumford regularity}.

\begin{thm}\label{T:CM-regularity}
For any coherent sheaf $\Gc\in\Coh(X)$ there exists an interger $m=m(\Gc)$ such that $H^j(X, \Gc(m))=0$ for all $j>0$ and all $d\geq m$. Moreover, the multiplication map $$
H^0(X,\Oc_X(i))\otimes H^0(X,\Gc(d))\longrightarrow H^0(X,\Gc(d+i))
$$
is surjective for all $i\geq 0$ and all $d\geq m$.
\end{thm}

More precisely, let $m\in \ZZ$. Recall that a coherent sheaf $\Gc$ on a polarized projective $\Spec(\mathbf{k})$-scheme $X$ is said to be \emph{$m$-regular} -- or of \emph{Castelnuovo-Mumford regularity} $m$ -- if
$$
H^{\alpha}(X,\Gc(m-\alpha))=0,
$$
for all integers $\alpha>0$. Now, if 
$$
0\rightarrow\Fc'\lra\Fc\lra\Fc''\rightarrow 0
$$
is a short exact sequence of coherent sheaves over $X$, additivity of the Euler characteristic on exact sequences implies that the regularity of $\Fc$ is bounded by the maximum of the regularities of $\Fc'$ and $\Fc''$ (see, e.g., \cite[Lecture 14]{Mum}).\\

Along the same lines, one can prove the following Theorem. We recall that the \emph{Hilbert polynomial} $h^{\Gc}$ of $\Gc\in \Coh(X)$ is defined as $h^{\Gc}(t):=\chi(\Gc(t))$. It is well known that $h^{\Gc}\in \QQ[t]$.

\begin{thm}[``Uniform Regularity Lemma'', \cite{Mum}]\label{T:uniform-regularity-lemma}
Let $\Gc$ be as above and let $h$ be a fixed Hilbert polynomial. Then the integer $m$ can be chosen so that all quotients of $\Gc$ with Hilbert polynomial $h$ and all of their kernels are $m$-regular.
\end{thm}

Finally, we collect in the following Theorem two fundamental results that we will use later on.

\begin{thm}\label{T:higher-direct-images-and-kuenneth}
(a) \cite{EGAIV2} Let $T$ be a $\mbfk$-scheme of finite type and let $\pi_T: X\times T\rightarrow T$ be the projection. If $\Fc$ is an $\Oc_T$-flat coherent sheaf on $X\times T$, then for $d\gg 0$ the direct images
$$
R^0{\pi_T}_*\Fc(d)
$$
are finite rank locally free sheaves on $T$.\\
(b) (see \cite{Mum} for the statement in this form) Let $W$ and $Z$ be any two algebraic $\mbfk$-schemes and let $\Hc$ and $\Kc$ be quasi-coherent sheaves on $W$ and $Z$, respectively. Denote by $\pi_W$ and $\pi_Z$ the respective projections from the product $W\times Z$, then 
$$
{\pi_W}_*(\pi_W^*\Hc\otimes\pi_Z^*\Kc)\simeq \Hc\otimes_{\mbfk} H^0(Z,\Kc).
$$
\end{thm}

\begin{rem}\label{R:grass-morphism-is-injective}
In the notation of Theorem \ref{T:higher-direct-images-and-kuenneth}, part (a), Theorem \ref{T:CM-regularity} implies that if $\Fc$ is a(n $m$-regular) quotient of ${\pi_X}^*\Gc$, then ${\pi_T}_*\Fc(m)$ uniquely determines the sheaf $\bigoplus_{d\geq m}{\pi_T}_*\Fc(d)$, which in turn determines $\Fc$ by the usual equivalence of categories between finitely generated graded modules and coherent sheaves on $X$ (\cite{FAC}).
\end{rem}

\section{A filtering cover for the schematic Grassmannian}\label{S:grass}

\subsection{Reminder on the classical construction}\label{SS:firstgrass}

We start by recalling the construction of the Grassmannian. For an integer $n\geq 1$, a scheme $S$ and quasi-coherent $\Oc_S$-module $\Ec$, we denote by $grass_n(\Ec)$ the set of locally free rank $n$ quotient $\Oc_S$-modules of $\Ec$.

\begin{thm}[\cite{EGAS}]\label{T:grass}
For every scheme $S$ and every quasi-coherent $\Oc_S$-module $\Ec$, the functor $\gamma_{n,\Ec}: (\Schb/S)^{\circ}\rightarrow \Setb$ given on objects by 
$$
\gamma_{n,\Ec} (T):= grass_n(\Ec_{T}),
$$
where $\Ec_{T}$ is the base change along the structure morphism $T\rightarrow S$, and as pullback on morphisms, is represented by a separated $S$-scheme $\Grass_n(\Ec)$. Moreover, there exists a locally free rank $n$ quotient $\Oc_{\Grass_n(\Ec)}$-module $\Qc$ of $\Ec_{\Grass_n(\Ec)}$, determined up to a unique isomorphism, such that 
$$
g\mapsto g^*(\Qc): \Hom_S(T,\Grass_n(\Ec))\buildrel{\sim}\over{\lra}\gamma_{n,\Ec}(T)
$$
is a natural isomorphism.
\end{thm}

The vector bundle $\Qc$ in the statement is the \emph{universal quotient bundle} of the Grassmannian. Note that if we do not assume that $\Ec$ in Theorem \ref{T:grass} is of finite type or finitely presented, then the scheme $\Grass_n(\Ec)$ will not in general be of finite type nor will it be quasi-compact. We will refer to $\Grass_n(\Ec)$ as the \emph{schematic}, or \emph{quasi-coherent}, Grassmannian when $\Ec$ is not assumed to have any finiteness properties. 
\vskip 3mm

Taking Theorem \ref{T:grass} for granted when $\Ec$ is a \emph{finitely presented} sheaf, we provide a construction of the schematic Grassmannian which is a filtering version of that of \cite{EGAS}. Our construction of the \emph{Quot} scheme of a quasi-coherent sheaf in the second part of this paper will partly follow the same pattern.\\

In the rest of this Section, we briefly review the part of the proof of Theorem \ref{T:grass} which we will need in the sequel. First, we have

\begin{lem}[\cite{EGAS}]
The functor $\gamma_{n,\Ec}$ is a sheaf of sets. 
\end{lem}

Therefore, we can 
reduce to proving its representability over the category of affine schemes. We will make such an assumption until the end of this section.

Let then $T$ be an $S$-scheme. For some index $i$, denote by $\gamma_{n,\Ec,i}(T)$ the subset of $\gamma_{n,\Ec}(T)$ consisting of the quotients $\Hc$ of $\Ec_T$ such that, for some finitely-presented subsheaf $\Ec^i$ of $\Ec$, we have a surjective composition

\begin{equation}\label{e:surjcompEiE}
\Ec^i_T\rightarrow \Ec_T \twoheadrightarrow \Hc,
\end{equation}
\noindent where the second arrow is the canonical quotient map. Note that the existence of an $i$ such that a surjective composition as in \eqref{e:surjcompEiE} exists follows from the Lemma below.

\begin{lem}{\cite[({\bf 0}, 5.2.3)]{EGAS}}
Let $T$ be a quasi-compact scheme and let $\Ec$ and $\Hc$ be two $\Oc_T$ modules. Assume $\Hc$ is of finite type and let moreover $u:\Ec\rightarrow \Hc$ be a surjective homomorphism. If $\Ec$ is a filtering inductive limit of a system $(\Ec^i)_i$ of $\Oc_T$-modules, then there exists an index $i$ such that there is a surjection $\Ec^i\rightarrow \Hc$.
\end{lem}


Therefore we have that
$$
T\mapsto\gamma_{n,\Ec,i}(T),
$$
together with the usual pull-back on morphisms, defines a subfunctor of $\gamma_{n,\Ec}$.

\begin{nota} 
To simplify the notation, from now on we will write $\gamma$ and $\gamma_i$ instead of $\gamma_{n,\Ec}$ and $\gamma_{n,\Ec,i}$, respectively, if no confusion can arise.
\end{nota}

Since we can now assume that $S$ is affine, the sheaf $\Ec$ is completely determined by a $\Gamma(S,\Oc_S)$-module $E$ via Serre's equivalence of categories: $\Ec=\tilde{E}$. Thus $\Ec$ is generated by a (possibly infinite) family of sections $(t_a)_{a\in \Omega}$. For an $S$-scheme $T$, denote by $t_{a,T}$ the pullback of $t_a$ along the structure morphism $T\rightarrow S$. Let then $H$ be a subset of $\Omega$ consisting of $n$ elements, using the sections $(t_{a,T})_{a\in H}$ we can define a homomorphism of $\Oc_T$-modules
$$
\varphi_{H,T}:\Oc_T^n\rightarrow \Ec_T.
$$
Now, consider the subset $F_H(T)$ of $\gamma(T)$ consisting of the quotients $\Hc$ of $\Ec_T$ such that we have a surjective composition

\begin{equation}\label{e:surjcompOE}
 \Oc_T^n\xrightarrow{\varphi_{H,T}} \Ec_T\stackrel{q}{\longrightarrow} \Hc,
\end{equation}

\noindent where the second arrow is the canonical quotient map. For future reference we denote by $s_a$ \label{sis} the canonical image in $\Gamma(T,\Hc)$ of the section $t_{a,T}$.\\
The datum $T\mapsto F_H(T)$ together with pullback on morphisms defines a subfunctor of $\gamma$, and the main step in Grothendieck's construction consists in proving that such a functor is represented by a scheme $X_H$ (which one could call the \emph{inverse Pl\"ucker subscheme}) which is affine over $S$, and that $F_H$ is an open subfunctor of $\gamma$. The functors $\Gc_i^j$ in the Lemma \ref{L:affmorphs} below are essentially unions of the functors $F_H$ as $H$ ranges over the set of sections of $\Ec^i$.\\

\subsection{The schematic Grassmannian as a filtering inductive limit}

The following Lemma in the spirit of \cite{EGAS} collects a few results which we will need in the rest of the paper.

\begin{lem}\label{L:surjptwise}
Let $Z$ and $Z'$ be two locally ringed spaces  and let $u:\Fc\rightarrow\Gc$ be a homomorphism of quasi-coherent $\Oc_Z$-modules of finite presentation. Then the following statements hold.
\item[(a)] \label{L:surjptwisea} The set of points $z$ of $Z$ where the localization $u_z:\Fc_z\rightarrow\Gc_z$ is surjective is open in $Z$.
\item[(b)] \label{L:surjptwiseb} The homomorphism $u_z:\Fc_z\rightarrow\Gc_z$ is surjective if and only if the homomorphism $u_z\otimes 1:\Fc_z/\mf{m}_z\Fc_z\rightarrow\Gc_z/\mf{m}_z\Gc_z$ is surjective.
\item[(c)] \label{L:surjptwisec} Let $f:Z'\rightarrow Z$ be a morphism of locally ringed spaces and put $\Fc'=f^*(\Fc)$, $\Gc'=f^*(\Gc)$ and $u'=f^*(u):\Fc'\rightarrow\Gc'$. Then the localization $u'_{z'}$ at a point $z'$ of $Z'$ is surjective if and only if the localization $u_z$ is surjective at the point $z=f(z')$.
\end{lem}

\noindent {\sl Proof:}
(a) Assume $u_z$ to be surjective at the point $z$. We will find a neighborhood $V$ of $z$ such that $u_{z''}$ is surjective at $z''$ for all $z''\in V$. For this, since $\Fc$ and $\Gc$ are of finite presentation, there exists a neighborhood $U$ of $z$ such that we have a commutative ladder diagram with exact rows
\[
          \begin{CD} 
          \Oc_Z^{s}|U         @>>>      \Oc_Z^{r}|U             @> v >> \Fc @>>> 0\\
          @VVV                                 @VVV                              @VV u V \\
          \Oc_Z^{s_1}|U  @>>>      \Oc_Z^{r_1}|U      @>>> \Gc @>>> 0,
          \end{CD}  \]
for some positive integers $r, r_1, s, s_1$.
By exactness of the top row, the surjectivity of $u$ is equivalent to the surjectivity of the composition $\varphi:=u\circ v: \Oc_Z^r|U\rightarrow \Gc$ and the same statement holds for the localized homomorphism $\varphi_z$. In order to conclude, we claim that $\varphi_z$ is surjective if and only if there exists a neighborhood $V$ of $z$ such that the localization $\varphi_{z''}$ is surjective for all $z''$ in $V$.\\
The question being local, consider instead a ring $A$, an $A$-module $N$, and a homomorphism $\varphi: A^r\rightarrow N$ such that at some point $\mf{p}$ the elements $n_1,\ldots,n_r\in N$ generate $N_{\mf{p}}$ over $A_{\mf{p}}$. Moreover, let $l_1,\ldots,l_t$ be generators for $N$ so that we have
\be\label{e:localz}
l_i=\sum_{j=1}^{r} a_{ij}n_j
\ee
where $a_{ij}\in A_{\mf{p}}$ for all $i$ and $j$.\\
Now, the localization of $A$ at $\mf{p}$ is a filtering inductive limit 
$$
A_{\mf{p}}=A[(A\setminus \mf{p})^{-1}]=\ind_{S\not\supset \mf{p}} A[S^{-1}]
$$
and, since the $a_{ij}$'s are finite in number, there exists a multiplicatively closed subset $S_0\not\supset \mf{p}$ of $A$ such that $a_{ij}\in A[S_0^{-1}]$ for all $i,j$. It follows that we can localize \eqref{e:localz} again at other points around $\mf{p}$ and get every time a surjective localized homomorphism.\\
(b) The question is again local, so consider a homomorphism of modules $u_{\mf{p}}: M_{\mf{p}}\rightarrow N_{\mf{p}}$ and the canonical quotient morphism $w_N:N_{\mf{p}}\rightarrow N_{\mf{p}}/\mf{p}N_{\mf{p}}$. Suppose $u_{\mf{p}}\otimes 1: M_{\mf{p}}/\mf{p}M_{\mf{p}}\rightarrow N_{\mf{p}}/\mf{p}N_{\mf{p}}$ is surjective. Then, composing the other canonical surjection $w_M:M_{\mf{p}}\rightarrow M_{\mf{p}}/\mf{p}M_{\mf{p}}$ with $u_{\mf{p}}\otimes 1$ we get a surjective homomorphism $M_{\mf{p}}\rightarrow N_{\mf{p}}/\mf{p}N_{\mf{p}}$. By Nakayama's Lemma, the module $N_{\mf{p}}$ has the same generators of $N_{\mf{p}}/\mf{p}N_{\mf{p}}$, whence our claim. The other direction of the argument is clear.\\
(c) We have $\Fc'_{z'}=\Fc_z\otimes_{\Oc_{Z,z}}\Oc_{Z,z'}$ and $\Gc'_{z'}=\Gc_z\otimes_{\Oc_{Z,z}}\Oc_{Z,z'}$ and $u'_{z'}$ is obtained from $u_{z}$ by base change from $\Oc_{Z,z}$ to $\Oc_{Z,z'}$. Let $\kappa$ and $\kappa'$ be the residue fields of $z\in Z$ and $z'\in Z'$, then $\Fc'_{z'}\otimes \kappa'=(\Fc_{z}\otimes \kappa')\otimes \kappa'$ and $\Gc'_{z'}\otimes \kappa'=(\Gc_{z}\otimes \kappa')\otimes \kappa'$, and we get $u_{z'}\otimes 1_{\kappa'}:\Fc'_{z'}\otimes \kappa'\rightarrow \Gc'_{z'}\otimes \kappa'$ by base changing $u_{z}\otimes 1_{\kappa}:\Fc_{z}\otimes \kappa\rightarrow \Gc_{z}\otimes \kappa$ from $\kappa$ to $\kappa'$. Now, since such a base change is faithfully flat one can conclude by applying Nakayama's Lemma and part (b).
\qed

\begin{lem}\label{L:affmorphs}
Let $S$ and $T$ be as in Section \ref{SS:firstgrass} and let $\Ec^i\buildrel{\alpha^{i,j}}\over{\longrightarrow}\Ec^j$ be a homomorphism of finitely presented $\Oc_S$-modules. Denote by $q:\Ec^j_T\longrightarrow\Hc$ the canonical quotient map. Then the functor 
$$
\Gc_i^j:(\Schb/S)^{\circ}\rightarrow \Setb
$$
defined by
\begin{equation}\label{e:surjcomp}
\Gc_i^j(T):=\gamma_{n,\Ec^j,i}(T)=\left\{ \Hc\in\gamma_{n,\Ec^j}(T)|\begin{array}{cc}{the}\,\,{composition}\,\,\,\Ec^i_T\buildrel{\alpha^{i,j}_T}\over{\longrightarrow} \Ec^j_T \buildrel{q}\over{\longrightarrow}\Hc\\
{is}\,\,{surjective}
\end{array}\right\}
\end{equation}

is an open subfunctor of $\gamma_{n,\Ec^j}$. Moreover, if we let $G_i^j$ be the open subscheme of $\Grass_n(\Ec^j)$ representing the above subfunctor, we have that there is an affine morphism
$$
v_{ij}: G_i^j\rightarrow\Grass_n(\Ec^i).
$$
\end{lem}

\noindent {\sl Proof:}
We start by proving that $\Gc_i^j$ is an open subfunctor of $\gamma_{n,\Ec^j}$. Let thus $Z$ be an $S$-scheme, we need to show that the fiber product functor
$$
T\mapsto\Gc_i^j(T)\times_{\gamma_{n,\Ec^j}(T)}\Hom_S(T,Z)
$$
is represented by an open subscheme of $Z$. Yoneda's lemma implies that a natural transformation $\Hom_S(-,Z)\Rightarrow\gamma_{n,\Ec^j}$ is completely determined by an element $\Fc\in\gamma_{n,\Ec^j}(Z)$ as a pullback: $\Hom_S(T,Z)\ni g\mapsto g^*(\Fc)\in \gamma_{n,\Ec^j}(T)$. By Lemma \ref{L:surjptwisea}(a) we have that the set of points of $Z$ where the localization of the composition $q\circ \alpha^{i,j}_Z$ is surjective is an open subset

\begin{equation}\label{e:U}
U:=U_{Z,i,\Fc}
\end{equation}

of $Z$. Moreover if $Y$ is another scheme, the set of $S$-morphisms $g:Y\rightarrow Z$ such that $g^*(\Fc)\in \Gc_i^j(Y)$ is equal to the set of $S$-morphisms such that $g(Y)\subset U$. In fact, the set of points $y$ of $Y$ where $(g^*(q)\circ \alpha^{i,j}_Z)_y$ is surjective is equal to $g^{-1}(U)$ by Lemma \ref{L:surjptwise}(c). On the other hand, saying that $g^*(\Fc)\in \Gc_i^j(Y)$ means that $g^{-1}(U)$ must coincide with all of $Y$.
We have just proved that the above fiber product of functors is represented by an open subscheme of $Z$.\\
\indent To establish the second part of the statement, consider the natural transformation
$$
\Gc^j_i\Longrightarrow\mathcal \gamma_{n,\Ec^i}
$$
defined by sending the quotients in $\Gc_i^j(T)$ to the corresponding elements of $\gamma_{n,\Ec^i}(T)$. We will show that such morphism of functors is represented by the affine morphisms of schemes $v_{ij}$ of the statement. To see this, it remains to show that the morphism is indeed affine. Now, replacing $\Ec$ with $\Ec^i$ in \eqref{e:surjcompOE}, we still obtain subfunctors $F_H^i$ of $\gamma_{n,\Ec^i}$ which are represented by affine subschemes $X_H^{i}$. Since we are assuming $S$ to be affine, the schemes $X_H^{i}$, as $H$ varies over the subsets of $(t_{a,T})_{a\in I}$ of cardinality $n$, can be identified with subschemes forming an open covering of $\Grass_n(\Ec^i)$. Given one of such schemes, which we denote again $X_H^{i}$ by a slight abuse of notation, we have to show that its inverse image is an open affine subscheme of $G_i^j$. At this point, given the homomorphism $\Ec^i\buildrel{\alpha^{i,j}}\over{\longrightarrow}\Ec^j$ over $S$ and a morphism $T\rightarrow S$, the fact that the composition $$
\Ec^i_T\buildrel{\alpha^{i,j}_T}\over {\longrightarrow}\Ec^j_T\br{q}\over\longrightarrow \Hc
$$
is surjective implies that the inverse image $v_{ij}^{-1}(X_H^{i})$ of $X_H^{i}\subset\Grass_n(\Ec^i)$ is equal to the subscheme $X_H^{j}$ of $G_i^j$. \qed\\

\begin{rem}
Note that the open subset $U_{Z,i,\Fc}$ in the above proof is equal to the union of the sets $U_{Z,H,\Fc}$ of \cite[Lemme {\bf 1},9.7.4.6]{EGAS}, as $H$ varies in the family of subsets of cardinality $n$ of $(t_{a,Z})_{a\in I}$. See also Example \ref{E:qcptplckr}.
\end{rem}

\begin{lem}\label{L:grsub}
The functor $\gamma_i$ is represented by the quasi-compact scheme given by the projective limit
$$
\Grqc(n,\Ec)_i:=\varprojlim_{j>i} G_i^j.
$$
\end{lem}

\noindent {\sl Proof:} First, we show that $\Grqc(n,\Ec)_i$ is quasi-compact.
The morphism in the filtering projective system 
\be\label{e:filtprojsysgrass}
\left((G_i^j)_{j\geq i}, G_i^j\buildrel{v_{j,j'}}\over\longleftarrow G_i^{j'}\right)
\ee
whose target is $G_i^i=\Grass_n(\Ec^i)$, is affine by Lemma \ref{L:affmorphs}. Moreover, all of the other morphisms in the system can be defined via natural transformations as in the proof of the previous Lemma, and proved to be affine arguing similarly. Therefore, quasi-compactness of the projective limit follows from Proposition \ref{P:qcam}.
\\

\indent In order to prove that $\gamma_{i}$ is representable, we have to show that the functors $\Hom_S(-,\Grqc(n,\Ec)_i)$ and $\gamma_{i}$ are naturally isomorphic. Using the isomorphism
$$
\Gc_i^j\tilde{\Longrightarrow} \Hom_S(-,G_i^j)
$$
and, for any $S$-scheme $T$,
$$
\Hom_S(T,\Grqc(n,\Ec)_i)=\Hom_S(T,\varprojlim_{j\geq i}G_i^j)=\varprojlim_{j\geq i}\vphantom{.}^{\Setb}\Hom_S(T,G_i^j),
$$
one can see that $\gamma_{i}(T)$ is the vertex of a left cone over the diagram formed by the sets $\left(\Hom_S(T,G_i^j)\right)_j$, with morphisms resulting from those of the system \eqref{e:filtprojsysgrass}. In this way, for every $T$ one gets a map 
$$
\gamma_i(T)\buildrel{L_T}\over{\longrightarrow}\varprojlim_{j\geq i}\Hom_S(T,G_i^j),
$$
which is natural in $T$ by the universal property of the projective limit.\\

Furthermore, we have that a right inverse $Y$ for $L$ is provided by Yoneda's Lemma, which implies that any natural morphism from $\Hom_S( - ,\Grqc(n,\Ec)_i)$ to $\gamma_{i}$ is completely determined by pulling back an element of $\gamma_{i}(\Grqc(n,\Ec)_i)$. Indeed, the square
\[
          \begin{CD} 
          \Hom_S(T',\varprojlim_{j\geq i}G_i^j)       @>Y_{T'}>>        \gamma_{i}(T')\\
          @VVV             @VVV \\
          \Hom_S(T,\varprojlim_{j\geq i}G_i^j)       @>Y_{T}>>  \gamma_{i}(T)
          \end{CD}  \]
is commutative, thanks to the fact that pulling back anticommutes with the composition of morphisms. At this point, one can see by direct computation that, for all $T$, the composition $L_T\circ Y_T$ of the two natural transformations is the identity.\\
\indent On the other hand, $Y$ is also a left inverse for $L$. To see this, let $\Fc^j$ be the universal sheaf over $G_i^j$, for $j\geq i$, and let $v_j:\Grqc(n,\Ec)_i\rightarrow G_i^j$ be the canonical maps. If $\Fc$ denotes the universal sheaf over $\Grqc(n,\Ec)_i$, we have
$$
\Fc=\varinjlim_{j\geq i} v_j^*\Fc^j.
$$
This allows us to conclude that the composition $Y_T\circ L_T$ is the identity.
\qed\\

\begin{lem}\label{L:opensubf}
 $\gamma_{i}$ is an open subfunctor of $\gamma$.
\end{lem}

\noindent {\sl Proof:} The same argument we used in the proof of Lemma \ref{L:affmorphs} to prove that $\Gc^j_i$ is an open subfunctor of $\gamma_{n,\Ec^j}$ can be used to show that $\gamma_{i}$ is an open subfunctor of $\gamma(-):=grass_n(\Ec_{(-)})$. 
\qed\\

\begin{lem}\label{L:grassopenemb}
 For $i<i'$ we have an open embedding of quasi-compact schemes
$$
\Grqc(n,\Ec)_i\rightarrow \Grqc(n,\Ec)_{i'}.
$$
\end{lem}

\noindent {\sl Proof:} From Lemma \ref{L:grsub} we know that for any $i$ the quasi-compact scheme $\Grqc(n,\Ec)_i$ represents the functor $\gamma_{i}$. Therefore, proving the statement amounts to showing that $\gamma_{i}$ is an open subfunctor of $\gamma_{i'}$ whenever $i<i'$, namely, that the functor $(\Schb/S)^{\circ}\rightarrow \Setb$ given by
$$
T\mapsto \gamma_{i}(T)\times_{\gamma_{i'}(T)}\Hom_S(T,Z),
$$
is represented by an open subscheme of $Z$. As we did in the previous proofs, after applying Yoneda's Lemma, the main tool we use is the following variation of \cite[Lemme \textbf{1},9.7.4.6]{EGAS} whose proof can be obtained in essentially the same way.

\begin{lem}\label{L:lemmevar}
\begin{itemize}
\item[(1)] Let $Z$ be an $S$-scheme, $\Fc$ a quotient $\Oc_Z$-module $\Ec_Z\buildrel{q}\over{\longrightarrow}\Fc$ of $\Ec_Z$ such that the composition
$$
\Ec_Z^{i'}\buildrel{\alpha^{i'}_Z}\over\longrightarrow \Ec_Z\buildrel{q}\over{\longrightarrow} \Fc
$$
is surjective. Then the set $U_{Z,i\rightarrow i', \Fc}$, of points of $Z$ where the localization of the composition 
$$
\Ec_Z^i\buildrel{\alpha^{i,i'}_Z}\over\longrightarrow\Ec_Z^{i'}\buildrel{\alpha^{i'}_Z}\over\longrightarrow \Ec_Z\buildrel{q}\over{\longrightarrow} \Fc
$$
is surjective, is open in $Z$.
\item[(2)] Let $Y$ be another $S$-scheme. Then the set of $S$-morphisms $g:Y\rightarrow Z$ such that $g^*(\Fc)\in \gamma_{i}(Y)$ is the set of $S$-morphisms such that $g(Y)\subset U_{Z,i\rightarrow i', \Fc}$.
\end{itemize}
\end {lem}

A direct application of Lemma \ref{L:lemmevar} concludes the proof of Lemma \ref{L:grassopenemb}.\\

\noindent {\sl Proof of Lemma \ref{L:lemmevar}:} (1). Follows immediately from Lemma \ref{L:surjptwise}(a).\\
(2). Note that $g^*(q):\Ec_Y\rightarrow g^*(\Fc)$ is again a quotient homomorphism and that $g^*(\alpha^{i,i'}_Z)=\alpha^{i,i'}_Y$. By Lemma \ref{L:surjptwise}(c), the set of points $y$ of $Y$ where the localization of $g^*(q)\circ \alpha^{i'}_Y\circ \alpha^{i,i'}_Y=g^*(q)\circ \alpha^{i}_Y$ is surjective is thus equal to $g^{-1}(U_{Z,i,\Fc})\subset Y$, where $U_{Z,i,\Fc}$ is the open subset defined in \eqref{e:U}. Now, since $g^*(\Fc)\in \gamma_i(Y)$, it follows that $g^{-1}(U_{Z,i,\Fc})$ must coincide with $Y$.
\qed\\
\qed

\begin{prop}\label{P:lastgrass}
Let $\Ec$ be a quasi-coherent sheaf over the scheme $S$. Then, as $i$ varies, the functors $\gamma_{i}$ form an open covering of the functor $\gamma$. Furthermore, we have
$$
\Grass_n(\Ec)=\varinjlim_i \Grqc(n,\Ec)_i.
$$
\end{prop}

\noindent {\sl Proof:} The fact that each of the $\gamma_{i}$'s is an open subfunctor of $\gamma$ was established in Lemma \ref{L:opensubf}. As in the proofs of the previous Lemmas, for an $S$-scheme $Z$ let $U_{Z,i,\Fc}$ be the open subscheme of $Z$ representing the usual fiber product functor
$$
T\mapsto \gamma_{i}(T)\times_{\gamma(T)}\Hom_S(T,Z).
$$
We show that the $U_{Z,i,\Fc}$'s cover $Z$ as $i$ varies. It is enough to show that the statement holds on points. Let then $\Fc\in \gamma(Z)$ and let $z$ be a point of $Z$. We have in particular a locally free rank $n$ sheaf $\Fc_z$ generated by the localization at $z$ of the $n$ sections $s_a$ (which were introduced on page \pageref{sis}) and an $n$-dimensional $\kappa(z)$-vector space $\Fc\otimes_{\Oc_{Z,z}}\kappa(z)$ with basis the $s_a(z)$'s. This said, since $\Fc$ is a quotient of finite type of the inductive limit $\Ec$, there must exist an index $i$ and a surjection 
$$
\Ec_Z^i\otimes_{\Oc_{Z,z}}\kappa(z)\rightarrow \Fc\otimes_{\Oc_{Z,z}}\kappa(z).
$$
Thus, by Lemma \ref{L:surjptwise}(b), we obtain a surjection $\Ec^i_{Z,z}\rightarrow \Fc_z$, hence $z\in U_{Z,i,\Fc}$ by definition.
This concludes our argument.
\qed\\

\begin{rem}
Note that Proposition \ref{P:lastgrass} implies that our construction of the Grassmannian is independent of the particular filtration of the sheaf $\Ec$ which we used.
\end{rem}

\begin{exa}\label{E:qcptplckr}
Here we will make use of the Pl\"ucker embedding, which is dealt with in detail in \cite[{\bf 1}, 9.8]{EGAS}.\\
When $S=\Spec({\mathbf k})$, for some field ${\mathbf k}$, we can give an explicit description of the schemes $\Grqc(n,\Ec)_i$ in terms of Pl\"ucker coordinates. For this, recall that in this case the sheaves $\Ec$ and $\Ec^i$ reduce to vector spaces, which we will denote by $E$ and $E^i$, respectively. Next, note that requiring the composition \eqref{e:surjcompEiE} to be surjective amounts to requiring that all of the compositions of the form
$$
\Oc_T^n\rightarrow E^i_T\rightarrow E_T \xrightarrowdbl{\mathsf {can}} \Hc
$$
be surjective, where the last arrow is the canonical surjection. Now, applying the Pl\"ucker functor we get that the corresponding compositions 
$$
\wedge^n\Oc_T^n\rightarrow \wedge^n E^i_T\rightarrow \wedge^n E_T {\twoheadrightarrow} \wedge^n\Hc,
$$
must be surjective, as well. That is, the scheme $\Grqc(n,\Ec)_i$ is determined by the non-vanishing of the $\binom{i}{n}$ Pl\"ucker coordinates whose multi-index contains only indices appearing among those of the basis vectors of $E_T^i$.
\end{exa}

\section{Representability of the quasi-coherent \emph{quot} functor.}

Throughout this Section, $S$ will be a noetherian scheme defined over a fixed algebraically closed field $\mbfk$, and $X$ will be a projective $S$-scheme (of finite type). By a coherent sheaf on $X$ we will mean a finitely presented quasi-coherent $\Oc_X$-module (see Lemma \ref{L:fpqconnoethiscoh}). 

\subsection{The coherent \emph{Quot} scheme}

Let now $T$ be another $S$-scheme, let $\pi_X: X\times_{S} T\rightarrow X$ be the projection and denote by $\Ec_T$ the pullback $\pi_X^*\Ec$, where $\Ec$ is a quasi-coherent $\Oc_X$-module. For a numerical polynomial $h\in \QQ[t]$, the $quot$ functor
$$
\eta_{h,\Ec}:= quot_h^X{\Ec_{(-)}}:({\Schb/S})^{\circ}\rightarrow{\Setb},
$$
is defined as
\begin{equation}\label{e:eta}
 \eta_{h,\Ec}(T)=\left\{\Kc\subset\Ec_T\,\, | \begin{array}{cc}\Ec_T/\Kc\,\, \textrm{is}\,\,\textrm{coherent,}\,\,\textrm{flat}\,\,\textrm{over}\,\,\Oc_T,\\\textrm{and}\,\,\textrm{has}\,\,\textrm{Hilbert}\,\,\textrm{polynomial}\,\, h
\end{array}\right\},
\end{equation}
together with pullback on morphisms. Without prescribing the Hilbert polynomial in the above definition one gets the general quot functor $\eta_{\Ec}$. 
Grothendieck's fundamental theorem reads as follows.

\begin{thm}[\cite{Gr}]\label{T:Grothmain}
Let $X$ be a projective $S$-scheme and let $\Gc$ be a coherent sheaf on $X$. Then, the functor $\eta_{h,\Gc}$ is represented by a projective $S$-scheme $\Quot_h(\Gc)$. Moreover, there exists a coherent quotient $\Qc\in \eta_{h,\Gc}(\Quot_h(\Gc))$ such that, for any $S$-scheme $T$, the morphism of functors
$$
\Hom_S(T,\Quot_h(\Gc))\ni g\longmapsto (\Id_X\times_S g)^*\Qc\in\eta_{h,\Gc}(T)
$$ 
is a natural isomorphism.
\end{thm}

We now briefly sketch the main idea in the proof of Theorem \ref{T:Grothmain}. From the results recalled in Section \ref{SS:prelquot}, we have that $m$-regularity of the coherent sheaf $\Gc\in\Coh(X)$ allows us to get, for any scheme $T$ over $S=\Spec(\mathbf k)$ and for any $T$-flat quotient homomorphism $\Gc_T\rightarrow \Fc$ with kernel $\Kc$, a short exact sequence of sheaves over $T$
$$
0\rightarrow \pi_{T*}\Kc(m)\lra H^0(X, \Gc(m))\otimes_{\mathbf k}\Oc_T\lra\pi_{T*}\Fc(m)\rightarrow 0,
$$
where $\pi_{T*}(\pi_{T}^*\Oc_T\otimes\pi_X^*\Gc(m))=H^0(X, \Gc(m))\otimes_{\mathbf k}\Oc_T$ by Theorem \ref{T:higher-direct-images-and-kuenneth}, part (b), and each of the sheaves in the above exact sequence is locally free by part (a) of the same Theorem.
By Remark \ref{R:grass-morphism-is-injective} we thus obtain an embedding of the functor $\eta_{h,\Gc}$ into the functor $\gamma_{h(m),H^0(X, \Gc(m))}$. The main step in the proof of Theorem \ref{T:Grothmain} is then the Theorem below (see \cite{Gr, Mum}).

\begin{thm}\label{T:stratum}
The scheme representing the functor $\eta_{h,\Gc}$ can be identified with a closed subscheme of the Grassmannian $\Grass_{h(m)} (H^0(X, \Gc(m)))$.
\end{thm}

In what follows we may sometimes refer to such a closed subscheme as the \emph{stratum} of the Grassmannian corresponding to the Hilbert polynomial $h$.

\vskip 5mm
Let now $\Ec\in\QCoh(X)$ be not necessarily coherent. We then have

\begin{lem}\label{L:sheafzt}
$\eta_{h,\Ec}$ is a sheaf in the Zariski topology on $\Schb/S$.
\end{lem}

\noindent {\sl Proof:} Let $\{U_{\alpha}\}_{\alpha}$ be a covering of the $S$-scheme $T$ and let $\Fc_{\alpha}\in \eta_{h,\Ec}(U_{\alpha})$. In the usual notation for restrictions, suppose that $\Fc_{\alpha,\beta}=\Fc_{\beta,\alpha}\in \eta_{h,\Ec}(U_{\alpha}\times_T U_{\beta})$, we want to find a unique sheaf $\Fc\in\eta_{h,\Ec}(T)$ whose restriction to $U_{\alpha}$ coincides with $\Fc_{\alpha}$.\\
For this, all we need to check is that if $U_{\alpha,\beta,\gamma}:=U_{\alpha}\times_T U_{\beta}\times_T U_{\gamma}$ is non empty, then the usual cocycle condition is satisfied. But this holds for the subsheaves, and hence for the quotients $\Fc_{U_{\alpha,\beta,\gamma}}\in \eta_{h,\Ec}(U_{\alpha}\times_T U_{\beta}\times_T U_{\gamma})$, since $\Ec_{T}$ itself is a sheaf. Moreover, $\Fc$ has Hilbert polynomial $h$ by semicontinuity, in particular by constance of the Hilbert polynomial on connected components, and flatness is a local condition.
\qed\\

Given a quasi-coherent sheaf $\Ec$ on $X$, our aim here is to construct an object $\Quot_h^X(\Ec)$, possibly in the category of $S$-schemes, that represents the functor $\eta_{h,\Ec}$.\\

Adapting pullbacks to the current $quot$ functor setting, in the notation of Lemma \ref{L:affmorphs} we can define the subfunctor

\begin{equation}\label{e:etai}
 \eta_{h,\Ec,i}(T):= \left\{ \Kc\in\eta_{h,\Ec}(T)\,\,|\begin{array}{cc}{the}\,\,{composition}\,\,\,\Ec_T^i\buildrel{\alpha^i_T}\over{\longrightarrow}\Ec_T\buildrel{q}\over{\rightarrow}\Ec_T/\Kc\\
{is}\,\,{surjective} \end{array}\right\},
\end{equation}

for an index $i$. Since $\eta_{h,\Ec,i}$ is a subfunctor of $\eta_{h,\Ec}$, it is also a sheaf of sets.

\subsection{Main results}

The following Lemma is the first main step in our construction.

\begin{lem}\label{L:qaffmorph}
Let $\Ec$ be a quasi-coherent $\Oc_X$-module.
Then for $i\leq a\leq b$ we have an affine morphism
$$
Q_i^b\rightarrow Q_i^a,
$$
from the scheme representing the functor $\eta_{h,\Ec^b,i}$ to the scheme representing $\eta_{h,\Ec^a,i}$. In particular, since $Q_i^i:=\Quot_h(\Ec^i)$, the morphism
$$
Q_i^a\rightarrow\Quot_h(\Ec^i)
$$
is affine.
\end{lem}

\noindent {\sl Proof:}
We have to show that the Grassmannian embedding of the $quot$ functor recalled in the previous section is compatible with our construction.\\
\indent First, note that if $\Ec^i$, $\Ec^a$ and $\Ec^b$ are three coherent sheaves on $X$, we can find a large enough integer $m$ such that all three of them are $m$-regular. Next, as we recalled above, $m$-regularity of the coherent sheaf $\Ec^i$ allows us to get, for any scheme $T$ over $S=\Spec(\mathbf k)$ and for any $T$-flat quotient homomorphism $\Ec^i_T\rightarrow \Fc$ with kernel $\Kc$, a short exact sequence of locally free sheaves over $T$
$$
0\rightarrow \pi_{T*}\Kc(m)\lra H^0(X, \Ec^i(m))\otimes_{\mathbf k}\Oc_T\lra\pi_{T*}\Fc(m)\rightarrow 0.
$$
We thus obtain an embedding of the functor $\eta_{h,\Ec^i}:=quot_h^X(\Ec^i_{(-)})$ into the functor $\gamma_{{h(m)},H^0(X,\Ec^i(m))}:=grass_{h(m)}{(H^0(X,\Ec^i(m))_{(-)})}$. This allows us to find a stratum of the Grassmannian that represents the functor $\eta_{h,\Ec^i}$ (Theorem \ref{T:stratum}).\\
Next, the homomorphism $\Ec^i\rightarrow \Ec^a$ induces a natural transformation
$$
\eta_{h,\Ec^a,i}\Rightarrow \eta_{h,\Ec^i},
$$
defined in the obvious way. The above transformation yields in turn a morphism of schemes
$$
Q_i^a\rightarrow Q_i^i,
$$
by representability of the \emph{quot} functor of a coherent sheaf and Yoneda's lemma. We claim that such a morphism is affine. We will use the covering induced on $\Quot_h(\Ec^i)$ by the one of the Grassmannian that was constructed in Section \ref{S:grass}.\\
In fact, thanks to Remark \ref{R:grass-morphism-is-injective}, from $\Ec^i\rightarrow \Ec^a$ we get a commutative square
\[
          \begin{CD} 
          Q_i^a       @>>>         Q_i^i \\
          @VVV             @VVV \\
          G_i^a       @>>>  G_i^i,
          \end{CD}  \]
where $G_i^a$ is the open part of the Grassmannian $\Grass_{h(m)}(H^0(X,\Ec^a(m)))$ whose points are isomorphism classes of quotients of $H^0(X,\Ec^i(m))$, and the vertical arrows denote the respective Grassmannian embeddings. By Lemma \ref{L:affmorphs} the lower arrow is an affine morphism, so we can conclude that $Q_i^a\rightarrow Q_i^i$ is also an affine morphism by restricting the lower arrow to the respective flattening stratum.\\
More generally, from the homomorphism $\Ec^a\rightarrow\Ec^b$, we get a natural transformation
$$
\eta_{h,\Ec^b,i}\Rightarrow \eta_{h,\Ec^a,i},
$$
and a resulting morphism of schemes $Q_i^b\rightarrow Q_i^a$. Keeping the notation as above we have another commutative diagram 
\[
          \begin{CD} 
          Q_i^b       @>>>         Q_i^a \\
          @VVV             @VVV \\
          G_i^b       @>>>  G_i^a,
          \end{CD}  \]
which allows us to conclude that $Q_i^b\rightarrow Q_i^a$ is affine, as well, by essentially the same argument.\\
Finally, let $S$ be any noetherian scheme over $\mbfk$. Then the statement follows from what we proved above plus the base change property of affine morphisms.
\qed\\

The next step in our construction is the Lemma below.

\begin{lem}\label{L:qqcpt}
Let $\Ec$ be a quasi-coherent $\Oc_X$-module. Then the functor $\eta_{h,\Ec,i}$ is represented by
$$
\Qqc(h,\Ec)_i:=\varprojlim_{a>i} Q_i^a,
$$
which is a quasi-compact scheme over $S$. 
\end{lem}

\noindent {\sl Proof:}
From Lemma \ref{L:qaffmorph} we see that all of the morphisms in the filtering projective diagram
\begin{equation}\label{e:amcd} 
\left((Q_i^a)_{a\geq i}, Q_i^a\leftarrow Q_i^{b}\right)
\end{equation}
are affine. As in the proof of Lemma \ref{L:grsub}, we then obtain that the projective limit is a quasi-compact scheme, by Proposition \ref{P:qcam}.\\
In order to conclude, it remains to prove that the scheme $\Qqc(h,\Ec)_i$ obtained as the projective limit of the diagram \eqref{e:amcd} actually represents the functor $\eta_{h,\Ec,i}$.
For this, the argument we used in Lemma \ref{L:grsub} for the functors $\gamma_{i}$ and the schemes $\Grqc(n,\Ec)_i$ still applies, provided that one uses Lemma \ref{L:qaffmorph} instead of Lemma \ref{L:affmorphs}.\qed
\vskip 1cm

\begin{lem}\label{L:openemb}
For $i\leq j$ we have an open embedding of schemes
$$
Q^a_i\rightarrow Q^a_j.
$$
\end{lem}

\noindent {\sl Proof:}
As usual, we prove the corresponding statement at the level of functors, i.e., we show that for every $S$-scheme $Z$ the fiber product functor 
$$
T\mapsto \eta_{h,\Ec^a,i}(T)\times_{\eta_{h,\Ec^a,j}(T)}\Hom_S(T,Z),
$$
is represented by an open subscheme of $Z$. Now, by definition of $\eta_{h,\Ec^a,i}$ we have a surjective composition
\begin{equation}\label{e:doublecomposite}
\Ec^i_Z\rightarrow\Ec^j_Z\rightarrow\Ec^a_Z\rightarrow\Fc,
\end{equation}
where the last homomorphism is the canonical quotient. 
Therefore the claim follows from Lemma \ref{L:lemmevar}, after applying the Grassmannian embedding to the composition \eqref{e:doublecomposite} for an $m$ large enough so that all the sheaves in question are $m$-regular.
\qed\\

Next, note that we have a commutative ladder diagram

\begin{equation}\label{e:ladderquot}
          \begin{CD}
          \vdots        @>>>  \vdots\\
          @VVV             @VVV \\
          Q_i^b      @>>>         Q_j^b \\
          @VVV             @VVV \\
          Q_i^a       @>>>  Q_j^a\\
          @VVV             @VVV \\
          \vdots        @>>>  \vdots\\
          \end{CD}          \\
\end{equation}

where the vertical arrows are surjections and the horizontal ones are open morphisms by Lemma \ref{L:openemb} above. Combining the Grassmannian embedding with the argument used in the proof of Lemma \ref{L:grassopenemb}, we obtain the following.

\begin{lem}\label{L:qopen}
Taking the projective limit over the upper indices in diagram \eqref{e:ladderquot} we get an open morphism of quasi-compact schemes
\be\label{e:qopen}
\Qqc(h,\Ec)_i \rightarrow \Qqc(h,\Ec)_j.
\ee
\end{lem}

Finally, define
\begin{equation}\label{e:qcQuot}
 \Quot_h(\Ec):=\varinjlim_i \Qqc(h,\Ec)_{i},
\end{equation}
where $\left(\Qqc(h,\Ec)_i\right)_i$ is the system of quasi-compact schemes and morphisms of the form \eqref{e:qopen} resulting from the above Lemmas.\\

\begin{thm} \label{T:iquot-ind-rep}
In the above notation, the functor $\eta_{h,\Ec}$ is covered by the functors $\eta_{h,\Ec,i}$.
\end{thm}

\noindent {\sl Proof:}
It remains to show that the subfunctors $\eta_{h,\Ec,i}$ cover $\eta_{h,\Ec}$ as $i$ varies. As in the case of the Grassmannian, it is enough to check this pointwise. Let $\Fc\in\eta_{h,\Ec}(Z)$, $z\in Z$, and consider the $\kappa(z)$-module of finite type $\Fc\otimes_{\Oc_{Z,z}}\kappa(z)$. Then there is an index $i$ such that we have a surjection 
$$
\Ec_Z^i\otimes_{\Oc_{Z,z}}\kappa(z)\rightarrow \Fc\otimes_{\Oc_{Z,z}}\kappa(z).
$$
At this point, the fact that there is a surjection $\Ec^i_{Z,z}\rightarrow \Fc_z$ follows from Nakayama's Lemma, since the quotient $\Fc$ is of finite type.
\qed\\

Now, the arguments that were used in Section \ref{S:grass} and what we have done so far in the current section yield, in combination with the Grassmannian embedding, that the functors $\eta_{h,\Ec,i}$ are open subfunctors of $\eta_{h,\Ec}$.
Taking Lemma \ref{L:qopen} into account we have the following.

\begin{thm} \label{T:iquot-rep}
Let $\Ec$ be a quasi-coherent sheaf over a projective $S$-scheme $X$. The functor $\eta_{h,\Ec}$ is represented by the scheme $\Quot^X_h(\Ec)$ from \eqref{e:qcQuot}.
\end{thm}

\noindent {\sl Proof:} As we said right before the statement, the functors $\eta_{h,\Ec,i}$ are open subfunctors of $\eta_{h,\Ec}$. Moreover, by Lemma \ref{L:qqcpt} such functors are representable and, by Lemma \ref{T:iquot-ind-rep}, they form an open covering of the functor $\eta_{h,\Ec}$. All of the above plus Lemma \ref{L:sheafzt} allow us to conclude.
\qed\\

The following Remarks and Example illustrate what we have achieved so far and relate the results to the literature.

\begin{rem}
From the above construction it follows that, when the sheaf $\Ec$ is not assumed to be coherent but just quasi-coherent, we obtain an infinite dimensional scheme $\Quot^X_h(\Ec)$ which, in principle, could be written as the $\Proj$ of some quasi-coherent algebra (see, e.g., \cite{EGAII}). Thus, even though infinite dimensional, our moduli space is an actual scheme and not an ind-scheme in the strict sense.
\end{rem}

\begin{exa}
Let $X=S$ in Theorem \ref{T:iquot-rep} above. Then, $\Quot_h^{X=S}(\Ec)$ reduces to a (relative) schematic Grassmannian of quotients of $\Ec$ of a certain rank prescribed by the Hilbert polynomial $h$ which, in this case, reduces to a constant.\\
In particular, let $\mathbf k$ be a field and let $X=S=\Spec({\mathbf k})$. An object $\Vc\in QCoh(X)$ is then a (possibly infinite dimensional) vector space over ${\mathbf k}$ and the scheme $\Quot_h^{\Spec({\mathbf k})}(\Vc)$ is then nothing but the usual Grassmannian $\Grass_h(\Vc)$.
\end{exa}

\begin{rem}\label{R:Kleiman}
More generally, let $S$ be an arbitrary scheme and $X$ be an $S$-scheme, not necessarily equal to $S$. Further, let $\Ec\in\QCoh(X)$. Then the definition of the \emph{quot} functor still makes sense.\\
In \cite[Prop. 2.2]{Kle} the author proves that if we consider the \emph{quot} functor of length $h\equiv 1$ quotients, we obtain that $\eta_{1,\Ec}$ is represented by the scheme
$$
\PP(\Ec):=\Grass_1(\Ec).
$$
That is, the $Quot$ scheme provides yet another way to define the projectivization $\PP(\Ec)\rightarrow X$ of a quasi-coherent sheaf.\\
It is worth mentioning that, for length $1$ quotients, Kleiman is able to show representability making essentially no assumption on $X$ and $S$, by exploiting an isomorphism between the $quot$ and Grassmann functors.
\end{rem}

\subsection{\emph{Uniformly regular} sheaves and a ``large scale''\\Grassmannian embedding}

Let again $\Ec$ be a quasi-coherent $\Oc_X$-module, and let $S=\Spec(\mbfk)$. We will show that in this case it is possible to obtain an analog of the classical Grassmannian embedding.\\

Motivated by the discussion in Section \ref{SS:prelquot}, we make the following definition.

\begin{Defi}\label{D:bcmg}
A quasi-coherent sheaf over a projective $\mathbf k$-scheme $X$ will be said to be \emph{uniformly m-regular} if there is an integer $m$ such that the Castelnuovo-Mumford regularities of its coherent approximations in the sense of Proposition \ref{P:Deligne} and Remark \ref{R:univincl} are all less than or equal to $m$.
\end{Defi}

\begin{lem}\label{L:seclast}
Let $\Ec$ be a uniformly $m$-regular quasi-coherent $\Oc_X$-module. Then there is a closed embedding
\be
\Qqc(h,\Ec)_i\rightarrow\Grqc\left(h(m),\ind_j H^0\left(X,\Ec^j(m)\right)\right)_i.
\ee
\end{lem}

\noindent {\sl Proof:} We go back to considering the components of the source and target schemes regarded as projective limits. In our usual notation, we have a commutative ladder diagram
\[
          \begin{CD} 
          \vdots        @>>>  \vdots\\
          @VVV             @VVV \\
          Q_i^b      @>>>         G_i^b \\
          @VVV             @VVV \\
          Q_i^a       @>>>  G_i^a\\
          @VVV             @VVV \\
          \vdots        @>>>  \vdots\\
          \end{CD}  \]\\
          
where the vertical arrows are surjective affine morphisms and the horizontal ones are the restrictions of the respective Grassmannian embeddings. The vertical morphisms being affine, we can reduce to proving the statement locally. \\
Let then $A$ and $B$ be two rings such that
$$
A=\ind_{\beta}A^{\beta}, \qquad B=\ind_{\beta}B^{\beta},
$$
and suppose $A^{\beta}\buildrel{\psi_{\beta}}\over\longleftarrow B^{\beta}$ is a surjective (quotient) homomorphism for all $\beta$, i.e., $A^{\beta}=B^{\beta}/I^{\beta}$ where $I^{\beta}=\ker(\psi_{\beta})$. We can then realize $A$ as a global quotient of $B$ modulo the ideal 
$$
\ind_{\beta}I^{\beta}.
$$
This allows us to establish the statement.
\qed\\

Writing $\Qqc(i):=\Qqc(h,\Ec)_i$ and $\Grqc(i):=\Grqc\left(h(m),\ind_j H^0\left(X,\Ec^j(m)\right)\right)_i$ to simplify the notation, we thus have another commutative ladder diagram

\[
          \begin{CD} 
          \vdots        @>>>  \vdots\\
          @VVV             @VVV \\
          \Qqc(i)      @>>>         \Grqc(i) \\
          @VVV             @VVV \\
          \Qqc(i')        @>>>  \Grqc(i')\\
          @VVV             @VVV \\
          \vdots        @>>>  \vdots\\
          \end{CD}  \]\\

\noindent where the horizontal arrows are the closed embeddings resulting from Lemma \ref{L:seclast} and the vertical ones are the open embeddings resulting from Lemma \ref{L:qopen} and Lemma \ref{L:grassopenemb}, respectively. In analogy with the notion of quasi-projectivity in finite dimensions, we make the following definition.

\begin{Defi}\label{D:quasiclosed}
We call \emph{quasi-closed} and embedding of schemes resulting from a limit of a ladder diagram like the above one.
\end{Defi}

We thus have the following result.

\begin{prop}\label{P:urgrassemb}
Let $\Ec$ be a uniformly $m$-regular quasi-coherent sheaf on a $\mbfk$-projective scheme $X$. Then there is a quasi-closed embedding 
$$
\Quot^X_h(\Ec)\hookrightarrow \Grass_{h(m)}\left(\ind_j H^0\left(X,\Ec^j(m)\right)\right).
$$
\end{prop}

\noindent {\sl Proof:} Follows from the argument preceding the statement. \qed

\vfill
--------------\\
\noindent{\sf Institut des Hautes \'Etudes Scientifiques}, {email:} \texttt{gennaro.dibrino@aya.yale.edu}
\end{document}